
\input amstex
\documentstyle{amsppt}
\topmatter
\title On Kodaira energy of polarized log varieties
\endtitle
\author Takao FUJITA \endauthor
\address {Department of Mathematics
\newline
Tokyo Institute of Technology
\newline
Oh-okayama, Meguro, Tokyo
\newline
152 Japan
\newline
e-mail:fujita\@math.titech.ac.jp}
\endaddress
\endtopmatter

\define\dnl{\newline\newline}
\define\lra{{\longrightarrow}}
\define\nl{\newline}

\define\Div{{\roman{Div}}}

\define\Inf{{\roman{Inf}}}
\define\Int{{\roman{Int}}}

\define\Lim{{\roman{Lim}}}

\define\Pic{{\roman{Pic}}}
\define\Sing{{\roman{Sing}}}

\define\WP{{\roman{WPic}}}


\define\dm{{\roman{dim}}}

\define\SF{\Cal F}
\define\SG{\Cal G}
\define\SL{\Cal L}

\define\SO{\Cal O}
\define\SP{\Cal P}

\define\BC{{\Bbb C}}
\define\BP{{\Bbb P}}
\define\BQ{{\Bbb Q}}
\define\BR{{\Bbb R}}
\define\BZ{{\Bbb Z}}

\define\Conj{{\bf Conjecture. }}
\define\Cor{{\bf Corollary. }}
\define\Def{{\it Definition. }}
\define\Fact{{\it Fact. }}
\define\Lem{{\bf Lemma. }}
\define\Prf{{\it Proof. }}

\define\Rmk{{\it Remark. }}
\define\Th{{\bf Theorem. }}

\define\ke{{\kappa\epsilon}}

\document
\noindent{\bf Introduction}

By a log variety we mean a pair $(V,B)$ consisting of a normal variety $V$ and
a $\BQ$-Weil divisor $B=\sum b_iB_i$ such that $0\le b_i\le 1$ for each $i$.
If it has only log terminal singularities (see (1.6) for the precise
definition), the (log) canonical $\BQ$-bundle $K(V,B)=K_V+B$ of $(V,B)$ is
well-defined.
Given further a big $\BQ$-bundle $L$ on $V$, the (log) Kodaira energy of
$(V,B,L)$ is defined by
$$\ke(V,B,L)=-\Inf\{t\in\BQ{\ \vert\ }\kappa(K(V,B)+tL)\ge0\}.$$
In this paper we are mainly interested in the case $\ke<0$, or equivalently,
$K(V,B)$ is not pseudo-effective.

According to the classification philosophy, at least when $B=0$, $V$ should
admit a Fano fibration structure in such cases.
In \S2, by using Log Minimal Model Program which is available in dimension
$\le3$ at present (cf. [{\bf Sho}], [{\bf Ko}]), we establish the existence of
such a fibration in the polarized situation.
Namely, under some reasonable assumptions, some birational transform
$(V',B',L')$ of $(V,B,L)$ admits a fibration $\Phi: V'\lra W$ onto a normal
variety $W$ with $\dm W<\dm V$ such that $K(V',B')-\ke(V,B,L)L'=\Phi^*A$ for
some ample $\BQ$-bundle $A$ on $W$.
Such a fibration is unique up to some birational equivalence and every general
fiber of $\Phi$ is Fano.
In particular we have $\ke(V,B,L)\in\BQ$, generalizing a result in [{\bf B}].

In \S3, we study the set of possible values of Kodaira energies for any fixed
$n=\dm V$, called the spectrum set.
If we consider the case $B=0$ only and $V$ is allowed to have only terminal
singularities, the spectrum set seems to have no negative limit point (cf.
(3.2)).
Using Kawamata's result on the boundedness of $\BQ$-Fano 3-folds, we prove the
above spectrum conjecture for $n\le3$ under the additional asumption that $V$
is $\BQ$-factorial.
On the other hand, if $B\ne0$ or $V$ is allowed to have log terminal
singularities, the spectrum set becomes more complicated and has many negative
limit points, even when $n=2$.
We just provide a few examples of this sort and present a conjecture.

Of course, Mori-Kawamata theory is of fundamental importance in our method.
In fact we need a slightly improved version than usual, so \S1 is devoted to
this purpose.
The key is the notion of {\it log ample} $\BQ$-bundle.

Our results in this paper form a philosophical background of the classification
theory of polarized (log) varieties by Kodaira energy.
See [{\bf BS}], [{\bf F2}], [{\bf F3}] for precise classification results.
\dnl
{\bf \S1. Preliminaries}

Here we review some results from Mori-Kawamata-Shokurov theory and fix notation
and terminology.
Basically we follow the notation in [{\bf F2}] and [{\bf KMM}].

(1.1) By a {\it variety} we mean an irreducible reduced complete algebraic
space of finite type over the complex number field $\BC$.
It is assumed to be projective almost always in this paper.
The group of invertible sheaves (resp. Cartier divisors) on a variety $V$ is
denoted by $\Pic(V)$ (resp. $\Div(V)$).
A $\BQ$-{\it bundle} (resp. $\BQ$-{\it divisor}) is an element of
$\Pic(V)\otimes\BQ$ (resp. $\Div(V)\otimes\BQ$).
For a $\BQ$-divisor $D$, the $\BQ$-bundle determined by $D$ is denoted by
$[D]$, or simply by $D$ when confusion is impossible or harmless.

A $\BQ$-{\it Weil divisor} on a normal variety $V$ is a $\BQ$-linear
combination of prime Weil divisors.
It is said to be $\BQ$-Cartier if some positive multiple of it is a Cartier
divisor.
The integral part of a Weil divisor $D$ will be denoted by $\Int(D)$.

The set $\WP(V)$ of reflexive sheaves on $V$ of rank one forms a group, where
the sum of $\SF$, $\SG\in\WP(V)$ is defined to be the double dual of
$\SF\otimes\SG$.
We have $\SO(D)\in\WP(V)$ for any Weil divisor $D$ on $V$, and
$\omega_V\in\WP(V)$ for the canonical sheaf $\omega_V$ of $V$.
An element of $\WP(V)\otimes\BQ$ is called a $\BQ$-{\it Weil sheaf}, or simply
a $\BQ$-{\it sheaf}.
It is said to be $\BQ$-invertible if it belongs to the subgroup
$\Pic(V)\otimes\BQ$.
We say that $V$ is $\BQ$-{\it factorial} if every $\BQ$-sheaf is
$\BQ$-invertible, or equivalently, every Weil divisor is $\BQ$-Cartier.
In such a case we sometimes say ``globally $\BQ$-factorial'', since $V$ may
have non-$\BQ$-factorial singularities.

(1.2) For any invertible sheaf $\SL$ and any positive integer $m$, the Iitaka
dimension $\kappa(\SL^{\otimes m})$ is equal to $\kappa(\SL)$.
Hence $\kappa(L)$ is well-defined in a natural way for any $\BQ$-bundle $L$.
$L$ is said to be {\it big} if $\kappa(L)=\dm V$.
Similarly, the ampleness, nefness etc. of a $\BQ$-bundle are well-defined.

If $L$ is nef, it is big if and only if $L^{\dm V}>0$, but this is not a good
definition of bigness in general.

Given a surjective morphism $f: V\lra S$, we define the notion of $f$-bigness,
$f$-ampleness and $f$-nefness in a natural way.
In particular, if $S$ is projective, a $\BQ$-bundle $L$ on $V$ is $f$-big
(resp. $f$-ample, $f$-nef) if and only if $L+f^*H$ is big (resp. ample, nef) on
$V$ for some $H\in\Pic(S)$.
In such a case we also say that $L$ is relatively big (resp. ample, nef) over
$S$.

(1.3) {\bf Kodaira's Lemma.} {\sl Let $f: V\lra S$ and $L$ be as above and
suppose that $f$ is projective.
Then $L$ is $f$-big if and only if $L-E$ is $f$-ample for some effective
$\BQ$-divisor $E$.}

The proof is easy and well-known.
{}From this we obtain the following

(1.4) \Cor {\sl Let $\Phi: V\lra W$ be a surjective morphism of varieties.
Let $L$ be a $\Phi$-big $\BQ$-bundle on $V$.
Then the restriction of $L$ to any general fiber of $\Phi$ is big.}

(1.5) For a variety $V$, let $\SP$ be the set of pairs $(M, D)$ consisting of a
normal birational model $M$ of $V$ and a prime Weil divisor $D$ on $M$.
Two such pairs $(M_1, D_1)$ and $(M_2, D_2)$ correspond to the same discrete
valuation of the function field of $V$ if and only if there is another pair
$(M', D')$ with birational morphisms $\pi_j: M'\lra M_j$ such that
$\pi_j(D')=D_j$.
This is an equivalence relation in $\SP$, and the equivalence class will be
called a {\it place} of $V$.

If a place $P$ is represented by $(M, D)$ with $\pi: M\lra V$, the subvariety
$Y=\pi(D)$ of $V$ is independent of the choice of the representative pair, and
will be called the {\it locus} of $P$ on $V$.
In such a case we say that $P$ lies over $Y$.

For any subvariety $X$ not contained in $\Sing(V)$, let $\nu: V'\lra V$ be the
normalization of the blow-up of $V$ along $X$, and let $E'$ be the proper
transform on $V'$ of the exceptional divisor over $X$.
The place represented by $(V', E')$ will be called the primary place over $X$.

(1.6) A {\it log variety} is a pair $(V, B)$ consisting of a normal variety $V$
and a $\BQ$-Weil divisor $B=\sum b_iB_i$ on $V$ such that $0<b_i\le 1$, where
$B_i$'s are the prime components of $B$ (possibly $B=0$).
If $b_i=1$, $B_i$ is called an {\it outer boundary component} of $(V, B)$.
A subvariety $Y$ of $V$ of codimension $r$ is called an {\it outer boundary
strata} if there are outer boundary components $D_1, \cdots, D_r$ such that $Y$
is an irreducible component of $D_1\cap\cdots\cap D_r$.

$(V, B)$ is said to be {\it log smooth} at a point $p$ on $V$ if $V$ is smooth
at $p$ and the support of $B$ has only normal crossing singularity at $p$.
It is said to be log smooth if it is so at every point on $V$.
Thus, a prime component of $B$ may have singularities.
A log desingularization of $(V, B)$ is a log smooth pair $(M, D)$ together with
a birational morphism $\pi: M\lra V$ such that $\pi_*D=B$.
Such a log desingularization exists by virtue of Hironaka's theory.

For a log variety $(V, B)$, the $\BQ$-sheaf $\omega_V+B$ is called the
canonical $\BQ$-sheaf of $(V, B)$ and is denoted by $\omega(V, B)$.
It is called a canonical $\BQ$-bundle and is denoted by $K(V, B)$ if it is
$\BQ$-invertible.
If so, a number $a_P=a_P(V, B)$ is defined for any place $P$ of $V$ as follows:
Take a representative pair $(M, D)$ of $P$ with a birational morphism $\pi:
M\lra V$, such that $M$ is smooth.
Let $B'$ be the proper transform of $B$ on $M$.
Then $K(M, B')-\pi^*K(V, B)=\sum \mu_iE_i$ for some $\pi$-exceptional prime
divisors $E_i$ on $M$.
We set $a_P=\mu_i$ if $D=E_i$ for some $i$, or $a_P=0$ otherwise.
Clearly this is independent of the choice of $(M, D)$.
In fact, in order to define $a_P$, it suffices to assume that $\omega(V, B)$ is
$\BQ$-invertible in a neighborhood of a general point of the locus of $P$ on
$V$.

If $a_P\ge -1$ at any place $P$, we say that $(V, B)$ has only {\it log
canonical} singularities.
If furthermore the equality holds only for places lying over an outer boundary
strata at any general point of which $(V, B)$ is log smooth, we say that $(V,
B)$ has only {\it log terminal} singularities (or simply $(V, B)$ is log
terminal), and such a pair $(V, B)$ will be called a log terminal variety.

\Rmk
The terminology ``log terminal'' is used in slightly different senses in other
papers, and there are several notions of ``log terminal'' singularities (see
e.g. [{\bf Ko}]).
They are equivalent to each other when $\Int(B)=0$, but there are delicate
differences when $\Int(B)\ne0$.
Our definition is perhaps the weakest one among those preserving important
properties, but details are omitted here, since we don't need other
definitions.

(1.7) Later we will use the following simple fact.

\Lem {\sl Let $(V,B)$ be a log terminal variety and let $E$ be an effective
$\BQ$-divisor on $V$.
Suppose that there is no outer boundary strata of $(V, B)$ contained in the
support of $E$.
Then $(V, B+\delta E)$ is log terminal for any sufficiently small $\delta>0$.}

The proof is easy and is left to the reader.

(1.8) \Fact {\sl Our log terminal singularities are rational.}

We just sketch the idea of the proof, since we don't need this fact in this
paper.
Suppose that $(V,B)$ is log terminal at $p\in V$.
For the sake of simplicity we assume that $V$ is affine.
Then $\SO_V$ is ample and the Weil sheaf $\SO_V(\Int(B))$ is spanned by global
sections.
So the Weil divisor $\Int(B)$ is linearly equivalent to an effective Weil
divisor $D$ such that no outer boundary strata of $(V,B)$ is contained in $D$.
Putting $B'=B+\delta(D-\Int(B))$ for some small $\delta>0$, we infer that
$(V,B')$ is log terminal and $K(V,B')=K(V,B)$.
Since $\Int(B')=0$, $V$ has only rational singularities as in [{\bf KMM}].

(1.9) Let $(V, B)$ be a log terminal variety and let $f: V\lra W$ be a
surjective morphism.
A $\BQ$-bundle $L$ on $V$ is said to be {\it log $f$-ample} on $(V, B)$ if
there is an effective $\BQ$-divisor $E$ such that $L$ is $f$-nef, $L-E$ is
$f$-ample and $(V, B+E)$ is log terminal.
If $W$ is a point, we just say ``log ample''.

Any log ample $\BQ$-bundle is nef and big.
The converse is true if $\Int(B)=0$, but not always true in general.
Indeed, a log ample bundle must be nef and big not only on $V$ but also on any
outer boundary strata; and even this condition is not enough for the log
ampleness.
But we have the following

(1.10) \Lem {\sl Let $(V, B)$ be a log terminal variety and let $f: V\lra W$ be
a surjective projective morphism.
Let $\infty(f)$ be the union of curves $C$ in $V$ with $f(C)$ being a point,
and suppose that there is no outer boundary strata contained in $\infty(f)$.
Then $f^*A$ is log ample on $(V, B)$ for any ample $\BQ$-bundle $A$ on $W$.}

\Prf We may assume that $A$ is an ample line bundle.
Take an ample line bundle $H$ on $V$.
For any outer boundary strata $Y$, $lA_V-H$ is generated by global sections at
a general point $y$ of $Y$ for $l\gg 0$, since $f$ is finite on a neighborhood
of $y$.
Therefore we have $E\in{\ \vert\ }lA_V-H{\ \vert\ }$ for some $l\gg 0$ such
that there is no outer boundary strata contained in the support of $E$.
By (1.7), we see easily that $f^*A=A_V$ is log ample.

(1.11) {\bf Vanishing Theorem.} {\sl Let $(V, B)$ be a log terminal variety,
let $f: V\lra W$ be a surjective morphism and let $L$ be a line bundle on $V$.
Suppose that $L-K(V, B)$ is log $f$-ample.
Then $R^jf_*\SO(L)=0$ for any $j>0$.}

\Prf
For the sake of brevity, we just say ``ample'' in the sense ``relatively ample
over $W$'' for the moment.
Take an effective $\BQ$-divisor $E$ such that $L-K(V, B)-E$ is ample and $(V,
B+E)$ is log terminal.
Replacing $B$ by $B+E$, we may assume that $A=L-K(V,B)$ is ample.
Then, for some $a\gg 0$, $aA$ is very ample and the Weil sheaf
$\SO_V(\Int(B))\otimes\SO(aA)$ is spanned by global sections.
Hence, as in (1.8), we have an effective Cartier divisor $H$, a small number
$\delta>0$ and an effective Weil divisor $D$ such that $\Int(B)+H$ is linearly
equivalent to $D$, $A-\delta H$ is ample, and $(V,B')$ is log terminal for
$B'=B+\delta(D-\Int(B))$.
Then $\Int(B')=0$ and $L-K(V,B')=A-\delta H$ is ample.
Applying [{\bf KMM; Th.1-2-5}] on $(V,B')$, we obtain the desired assertion.

This generalizes the famous result of Kawamata-Viehweg.
Similarly we can generalize many results in [{\bf KMM}], replacing the
ampleness assumption by log ampleness.
For example

(1.12) {\bf Fibration Theorem} (compare [{\bf KMM; 3-1-1} \& {\bf 3-2-1}]).
{\sl Let $(V, B)$ be a log terminal variety and let $f: V\lra S$ be a
surjective morphism.
Let $L$ be an $f$-nef line bundle such that $mL-K(V, B)$ is log $f$-ample for
some $m>0$.
Then there is an $S$-morphism $\Phi: V\lra W$ onto a normal variety $W$ with
$g: W\lra S$ and a $g$-ample line bundle $A$ on $W$ such that $L=\Phi^*A$ and
$\Phi_*\SO_V=\SO_W$.}

\Prf
$mL-K(V, B)-E$ is $f$-ample and $(V, B+E)$ is log terminal for some $E$.
So we replace $(V, B)$ by $(V, B+E)$, and argue as usual.

(1.13) {\bf Rationality Theorem} (compare [{\bf KMM; 4-1-1}]).
{\sl Let $(V, B)$, $f: V\lra S$ be as above and let $H$ be a log $f$-ample line
bundle on $V$.
Then $\tau=\Inf\{t\ge0{\ \vert\ }K(V, B)+tH{\text{\ is
$f$-nef}}\}\in\BQ\cup\{\infty\}$.}

\Prf
$H-E$ is $f$-ample and $(V, B+E)$ is log terminal for some $E$.
We may assume that $K(V,B)+aH$ is not $f$-nef for some small $a>0$.
Since $K(V, B)+aH=K(V, B+aE)+a(H-E)$, there are only finitely many extremal
rays $R_j$ such that $(K(V,B)+aH)R_j<0$ by the Cone Theorem [{\bf KMM; 4-2-1}].
Then $K(V, B)+tH$ is $f$-nef if and only if $(K(V, B)+tH)R_j\ge0$ for each $j$,
hence $\tau\in\BQ\cup\{\infty\}$.

(1.14) In the sequel we shall freely use Log Minimal Model Program as in [{\bf
KMM}], [{\bf Ko}].
Terminologies such as extremal ray, contraction of it, divisorial contraction,
log flip and so on are used in the usual way.
Technical details are not necessary here, so we omit it.\dnl
{\bf \S2. Adjoint fibration}

(2.1) \Def
Let $L$ be a big $\BQ$-bundle on a log terminal variety $(V,B)$.
The {\it log Kodaira energy} of such a triple $(V,B,L)$ is defined as follows:
$$\ke(V,B,L)=-\Inf\{t\in\BQ{\ \vert\ }\kappa(K(V,B)+tL)=\dm V\}.$$

(2.2) The purpose of this section is to prove the following

\Th {\sl Let $(V,B)$ be a log terminal variety.
Suppose that $K(V,B)$ is not nef, $V$ is $\BQ$-factorial and $n=\dm V=3$.
Let $L$ be a big $\BQ$-bundle on $V$ such that $K(V,B)+aL$ is log ample on
$(V,B)$ for some $a>0$.
Then there is a birational transform $(V', B')$ of $(V,B)$ together with a
fibration $\Phi: V'\lra W$ such that\nl
1) $W$ is a normal variety with $\dm W<n$.\nl
2) $(V',B')$ is log terminal and $K(V',B')-\kappa\epsilon(V,B,L) L'=\Phi^*H$
for some ample $\BQ$-bundle $H$ on $W$, where $L'$ is the proper transform of
$L$ on $V'$ as a $\BQ$-Weil sheaf (corresponding to the proper transform of
Weil divisors).

In particular $\ke\in\BQ$.}

The proof consists of several steps.

(2.3) The transformation from $V$ to $V'$ is a sequence of elementary
divisorial contractions and log flips, as described below.

To begin with, set $\tau=\Inf\{t{\ \vert\ }K(V,B)+tL {\text{ is nef}}\}$.
Take $a$ such that $A=K(V,B)+aL$ is log ample.
By (1.13), $s=\Inf\{t{\ \vert\ }K(V,B)+tA {\text{ is nef}}\}\in\BQ$.
Hence $\tau=as/(1+s)\in\BQ$.

Applying (1.12) we get a fibration $f:V\lra X$ such that $K(V,B)+\tau L=f^*A$
for some ample $\BQ$-bundle $A$ on $X$.

(2.4) Now we let the Log Minimal Model Program run in this relative situation
over $X$:
If $K(V,B)$ is not $f$-nef, there is an extremal ray $R$ such that $K(V,B)R<0$
and $f_*R=0$.
Let $\rho: V\lra V^\flat$ be the contraction morphism of $R$, which is an
$X$-morphism.
If $\rho$ is of fibration type, then the Program ends.
If $\rho$ is a birational divisorial contraction, then $(V^\flat,B^\flat)$ is
log terminal, where $B^\flat=\rho_*B$.
Moreover $V^\flat$ is $\BQ$-factorial.
If $K(V^\flat,B^\flat)$ is nef over $X$, the Program ends.
Otherwise, replacing $(V,B)$ by $(V^\flat,B^\flat)$, we repeat the same
process.
If $\rho$ is a small contraction, then we take a log flip $V^+\lra V^\flat$ of
$\rho$ (cf. [{\bf Ko}], [{\bf Sho}]) and let $B^+$ be the proper transform of
$B$.
Then $(V^+, B^+)$ is log terminal and $V^+$ is $\BQ$-factorial.
If $K(V^+,B^+)$ is nef over $X$, the Program ends.
Otherwise we repeat the same process replacing $(V,B)$ by $(V^+,B^+)$.

As usual, by the termination theorem of log flips and by the decreasing
property of the Picard number, the Program must end after finite steps.
There are two possibilities for the final state:
Either we get a contraction $\rho$ of fibration type, or we get a model
$(V_1,B_1)$ over $X$ such that $K(V_1,B_1)$ is relatively nef over $X$.

(2.5) We will examine how the $\BQ$-sheaf $L$ (this may be viewed as if a Weil
divisor, if you like) during the above process (2.4).

Suppose that $\rho: V\lra V^\flat$ is a divisorial contraction.
Let $E$ be the exceptional divisor of $\rho$.
In this case $L^\flat=\rho_*L$ is a $\BQ$-bundle on $V^\flat$ such that
$\rho^*L^\flat=L+aE$ for some $a\in\BQ$.
We have $LR>0$ since $(K(V,B)+\tau L)R=0$.
Therefore $a>0$ since $ER<0$.
In particular $L^\flat$ is big as well as $L$, and $K(V^\flat,B^\flat)+\tau
L^\flat$ is the pull-back of $A$.
However, $L^\flat$ may not be nef even if $L$ is ample, and the invertibility
is not always preserved either.

Suppose that $\rho: V\lra V^\flat$ is a small contraction.
Let $L^+$ be the proper transform of $L$ on $V^+$.
Then $K(V^+,B^+)+\tau L^+$ is the proper transform of $K(V,B)+\tau L=A_V$, so
it is the pull-back of $A$ on $V^+$.
We claim that $L^+$ is big as well as $L$.

Indeed, for any place $P$, we have $a_P(V,B)\le a_P(V^+,B^+)$ (cf. [{\bf
KMM;5-1-11},(3)]), where $a_P$ is as in (1.6).
Hence $v_P(L^+-L)=\tau^{-1}v_P(K(V,B)-K(V^+,B^+))\ge0$, where $v_P$ is the
discrete valuation at $P$.
Thus $L^+-L$ is effective and $L^+$ is big.
Note that the strict inequality holds in the above situation if and only if the
locus of $P$ is contained in the exceptional set of $\rho$ (cf. [{\bf KMM;
5-1-11}]).

Thus, in either birational contractions, the bigness of $L$ is preserved.

(2.6) The Kodaira energy does not vary during the above process.

To see this, consider first the case of flip.
Since $K(V,B)+\tau L$ is nef, we have $\ke(V,B,L)\ge-\tau$.
Similarly $\ke(V^+,B^+,L^+)\ge-\tau$.
On some model $M$ we have an effective $\BQ$-divisor $E=L^+_M-L_M$.
For any $t\le\tau$, we have $K(V,B)+tL=K(V^+,B^+)+tL^++(\tau-t)E$, so
$\kappa(K(V,B)+tL)=\kappa(K(V^+,B^+)+tL^+)$ since $E$ is exceptional with
respect to $M\lra V^+$.
This implies $\ke(V^+,B^+,L^+)=\ke(V,B,L)$.

When $\rho$ is divisorial, the argument is similar and simpler.

(2.7) Now we consider the case in which we get a contraction $\rho: V'\lra
V^\flat$ of fibration type by the process (2.4).
Let $B'$ and $L'$ be the proper transforms of $B$ and $L$ on $V'$.
Then $K(V',B')+\tau L'=A_{V'}$, so $\kappa(K'+\tau L')=\dm X\le\dm V^\flat<n$.
Hence $\ke(V',B',L')=-\tau$.

By virtue of (2.6), the assertion of the Theorem (2.2) is satisfied for $W=X$.

(2.8) Next we consider the case in which we get a model $(V_1, B_1)$ such that
$K_1=K(V_1,B_1)$ is nef over $X$.
First we claim $\dm X=n$.

Indeed, otherwise, $L_1$ is big on any general fiber $F$ of $f_1: V_1\lra X$ by
(2.5) and (1.4).
Since $K_1$ is nef on $F$, $(K_1+\tau L_1)_F$ is big, but $K_1+\tau L_1=0$ on
$F$.
This is impossible unless $\dm F=0$, i.e., $\dm X=n$.

Next we claim that $K_1+\tau L_1=A_{V_1}$ is log ample on $(V_1, B_1)$.
To show this, we may assume that $L$ is log ample on $(V,B)$, replacing $L$ by
$K(V,B)+aL$ if necessary.
We will apply (1.10) for this purpose.
Let $Y$ be an outer boundary strata of $(V_1,B_1)$.
It suffices to derive a contradiction assuming $Y\subset\infty(f_1)$, where
$\infty(f_1)$ is the union of curves in $V_1$ which are contracted to points by
$f_1: V_1\lra X$.

Let $P$ be the primary place over $Y$.
If the locus $Y^\sharp$ of $P$ on $V$ is contained in the exceptional set of
the birational map $V\lra V_1$, then $a_P(V_1,B_1)>a_P(V,B)$ by the observation
(2.5).
But then $(V,B)$ is not log canonical since $a_P(V_1,B_1)=-1$, contradiction.
Therefore $Y^\sharp$ is not in the exceptional set, and $V\lra V_1$ is an
isomorphism at any general point of $Y^\sharp$.
In particular $Y^\sharp$ itself is an outer boundary strata of $(V,B)$.
Recalling $Y\subset\infty(f_1)$, we take a curve $C$ passing a general point of
$Y$ such that $f_1(C)$ is a point.
Let $C^\sharp$ be the proper transform of $C$ on $V$.
Now, take an effective $\BQ$-divisor $\Delta$ on $V$ such that $L-\Delta$ is
ample and $(V,B+\Delta)$ is log terminal.
Then $Y^\sharp\not\subset\Delta$ since otherwise $(V,B+\Delta)$ would not be
log terminal.
Hence $C^\sharp\not\subset\Delta$ and $LC^\sharp>\Delta C^\sharp\ge0$.
On the other hand, as we saw in (2.5), $L_1-L=\tau^{-1}(K(V,B)-K(V_1,B_1))$ is
an effective $\BQ$-divisor on some common model $\tilde V$ of $V$ and $V_1$,
whose locus on $V$ is contained in the exceptional set of $V\lra V_1$.
Hence $(L_1-L)\tilde C\ge 0$ for the common proper transform $\tilde C$ of $C$
and $C^\sharp$ on $\tilde V$.
Since $L_1\tilde C=L_1C$ and $L\tilde C=LC^\sharp$, we obtain $L_1C>0$.
But $f_1(C)$ is a point, so $K_1C\ge 0$ and $(K_1+\tau L_1)C=AC=0$.
Thus we get a contradiction, as desired.

\Rmk $f^*A$ is not always log ample on $(V,B)$.

(2.9) Since $A$ is ample on $X$, $K_1+lA$ is nef on $V_1$ for $l\gg0$.
So $\tau_2=\Inf\{t{\ \vert\ }K(V_1,B_1)+tL_1 {\text{ is nef}}\}<\tau$.
Since $K_1+{\tau}L_1=f_1^*A$ is log ample, we have $\tau_2\in\BQ$ and we get a
fibration $f_2: V_1\lra X_2$ such that $K_1+\tau_2L_1=f_2^*A_2$ for some ample
$\BQ$-bundle $A_2$ on $X_2$.
Now the situation is as in (2.3), and we let run the Log Minimal Model Program
run over $X_2$ as in (2.4).

If $\dm X_2<n$, we are done as in (2.7).
If $\dm X_2=n$, we get a model $(V_2,B_2,L_2)$ as before such that
$K_2=K(V_2,B_2)$ is nef over $X_2$.
$A_2$ is log ample on $(V_2,B_2)$, so we can repeat again.

We repeat this process as long as possible.
But the birational transformations $V\lra V_1\lra$ $ V_2\lra\cdots$ consist of
flips and divisorial contractions.
By the termination theorem of log flips and by the decreasing property of the
Picard number, as usual we infer that the whole process must terminate after
finite steps.
At the end we get a model with the desired fibration structure.
\comment
(2.10) \Cor {\sl Let $(V,B,L)$ be as in the theorem (2.2).
Then there is a log terminal variety $(F,D)$ and an ample $\BQ$-bundle $H$ on
$F$ such that $K(F,D)=kH$ for $k=\ke(V,B,L)$.}

\Prf
The final state of the preceding process is as in (2.7).
Namely, let $(V',B',L')$ be the triple as in the theorem (2.2) with $\Phi:
V'\lra W$.
We further let run the Program over $W$.
Then, after some preparatory birational transformations, we get a contraction
$\rho: V''\lra V^\flat$ of fibration type.
Let $F$ be a general fiber of $\rho$, $D$ be the restriction of the proper
transform of $B'$ and let $H$ be the restriction of the proper transform of
$L'$.
Since $\rho$ is the contraction of an extremal ray $R$, the bigness of $L'$
implies $L'R>0$, hence $L'$ is $\rho$-ample and $H$ is ample on $F$.
Thus $(F,D,H)$ has the desired property.
\endcomment
\dnl
{\bf \S3. Spectrum Conjecture}

(3.1) A {\it polarized terminal variety} is a pair $(V,L)$ consisting of a
variety $V$ having only terminal singularities and an ample line bundle $L$ on
$V$.
Its Kodaira energy is defined by
$$\ke(V,L)=-\Inf\{t\in\BQ{\ \vert\ }\kappa(K+tL)\ge0\}$$
as in (2.1), where $K$ is the canonical $\BQ$-bundle of $V$.
The spectrum set $S_n$ of polarized terminal $n$-folds is defined to be the set
of all the possible values of Kodaira energies of polarized terminal varieties
of dimension $n$.
It is easy to see $S_1=\{t\in\BQ{\ \vert\
}t\ge0\}\cup\{\frac{-2}d\}_{d\in\BZ^+}$.
Moreover, for any $n$, every positive rational number is contained in $S_n$.
As for the negative range, we have the following

(3.2) {\bf Spectrum Conjecture.} {\sl For any $\epsilon>0$, $\{t\in S_n{\
\vert\ }t<-\epsilon \}$ is a finite subset of $\BQ$.}

This is closely related to the following conjecture on the boundedness of
terminal Fano $n$-folds.

(3.3) \Conj {\sl For any fixed $n\in\BZ^+$, there exist positive integers $r$
and $d$ with the following properties: $(-K_V)^n<d$ and $rK_V$ is invertible
for every Fano $n$-fold $V$ which has only terminal singularities and Picard
number one.}

This is proved for $n=3$ by Kawamata [{\bf Ka}] under the additional assumption
that $V$ is $\BQ$-factorial.

(3.4) We will prove the Spectrum Conjecture for $n\le3$ using Kawamata's
result, at least for a $\BQ$-factorial variety $V$.
Let $S_n^o$ be the set of possible values of $\ke(V,L)$, where $V$ is
$\BQ$-factorial and has only terminal singularities, the Picard number of $V$
is one, and $L$ is a Weil sheaf which is ample as a $\BQ$-bundle.
Then

(3.4.1) {\sl $\{t\in S_n^o{\ \vert\ }t<-\epsilon\}$ is a finite subset of
$\BQ$.}

To see this, let $r$ and $d$ be as in (3.3).
Then for any $(V,L)$ as above, we have
$\ke(V,L)=-r^{n-1}(-K_V)^n/L(-rK_V)^{n-1}$, while $r^{n-1}(-K_V)^n<r^{n-1}d$
and $L(-rK_V)^{n-1}$ is a positive integer.
In the range $\ke<-\epsilon$,
$L(-rK_V)^{n-1}<\epsilon^{-1}r^{n-1}(-K_V)^n<\epsilon^{-1}r^{n-1}d$, so there
are at most finitely many such numbers.

(3.4.2) {\sl $\{t\in S_n{\ \vert\ }t<0\}\subset S_n^o\cup(\cup_{j<n}S_j)$, if
the Minimal Model Program works in dimension $n$.}

To see this, we apply (2.2) to the triple $(V,0,L)$.
Let $\Phi: V'\lra W$ be the fibration as there and let $L'$ be the proper
transform of $L$ on $V'$.
As in (2.7), we further let the Minimal Model Program run over $W$.
Then we get a contraction $\rho: V''\lra V^\flat$ of an extremal ray $R$ of
fibration type.
Let $F$ be a general fiber of $\rho$ and let $L''$ be the proper transform of
$L'$ on $V''$.
Then $K_F=\ke(V,L)L''_F$.
Moreover, $L''_F$ is big by (1.4), so $L''R>0$, which implies that $L''_F$ is
ample.
Hence $\ke(V,L)=\ke(F,L''_F)\in S_j$ for $j=\dm F$ if $j<n$.

If $j=n$, $V^\flat$ is a point and $V''=F$.
Moreover the Picard number of $V''$ is one.
Therefore $\ke(V,L)\in S_n^o$ in this case since $V''$ is still $\BQ$-factorial
and has only terminal singularities by the minimal model theory.
Thus we prove the claim (3.4.2).

Now, combining these observations (3.4.1) and (3.4.2), we obtain (3.2) by
induction on $n$.

It is uncertain whether this argument works for $n=3$ without the
$\BQ$-factoriality assumtion, but there is no such trouble in the cases $n<3$.

(3.5) \Rmk For $n>3$, we can show that $\{t\in S_n^s{\ \vert\ }t<3-n\}$ is a
finite set, where $S_n^s$ is the set of all the possible values of Kodaira
energies of polarized $n$-folds $(V,L)$ such that $V$ is smooth (cf. [{\bf
BS}], [{\bf F2}]).

(3.6) How about the spectrum set of polarized log varieties ?
The problem seems to be very complicated even in the following simplest cases.

Let $S_n^{ls}=\{\ke(V,B,L){\ \vert\ }(V,B)$ is a log smooth log variety,
$B=\Int(B)$, $L$ is an ample line bundle on $V\}$.
For $n=1$, this set is discrete in the negative range, but this is no more true
for $n\ge2$.

To see this, consider the case $V\cong\Sigma_e$ with $e>0$, the $e$-th
Hirzebruch surface, $B$ is the negative section with $B^2=-e$ and $L=B+(e+1)F$,
where $F$ is a fiber of the $\BP^1$-bundle map $\Sigma_e\lra \BP^1$.
Then the canonical bundle $K$ of $V$ is $-2B-(e+2)F$, hence
$K+B+(1+\frac1{e+1})L=(e+1)^{-1}B$ and $\ke(V,B,L)=-(1+\frac1{e+1})$.
Thus $S_2^{ls}$ contains $\{-1-\frac1{e+1}\}$, a sequence converging to $-1$
from below.
Moreover, we have $\ke(V,B,mL)=m^{-1}\ke(V,B,L)$ for any positive integer $m$,
so $\frac{-1}m$ is also a limit point of $S_2^{ls}$.
Thus $S_2^{ls}$ has infinitely many limit points in the negative range.

As in (3.4.2), we can show that $S_n^{ls}$ is contained in the union of low
dimensional log Kodaira spectrum sets and the set of values of Kodaira energies
of $(V,B,L)$, where $(V,B)$ is log terminal and $\Int(B)=B$, $V$ is
$\BQ$-factorial and of Picard number one, and $L$ is a Weil sheaf which is
ample as a $\BQ$-bundle.
However, unlike in (3.4), there is no counter part of (3.3), so the log
spectrum set has a complicated structure.
Here we just present the following conjecture, which was originally proposed by
Shokurov in a slightly different form:

(3.7) {\bf Log Spectrum Conjecture. } {\sl Given a subset $X$ of $\BR$, let
$\Lim(X)$ be the set of limit points of $X$, and put
$\Lim^k(X)=\Lim(\Lim^{k-1}(X))$.
Then $\Lim^k(S_n^{ls})\subset\{t\in\BR{\ \vert\ }t\ge k-n\}$ for any positive
integer $k\le n$.}

In particular, $\Lim(S_n^{ls})$ seems to resemble the set $S_{n-1}^{ls}$.
We conjecture also that the limit point can be reached only from below.
Namely, for any $x<0$, $\{t\in S_n^{ls}{\ \vert\ }x<t<x+\epsilon\}$ is a finite
set for some small $\epsilon>0$.

Similar phenomenon occurs if we allow $V$ to have log terminal singularities
instead of boundaries.
Moreover, the spectrum set of log terminal $n$-folds with $B=0$ seems to
resemble the above set $S_n^{ls}$, and to be closely connected with the set of
possible indices of log terminal $\BQ$-Fano $n$-folds.
See [{\bf A}].

(3.8) As for the lowest possible value of $\ke$, we have the following fact:\nl
{\sl Let $L$ be an ample line bundle on a log terminal $n$-fold $(V,B)$.
Then $\ke(V,B,L)\ge-n$ unless $V\cong\BP^n$.}

In fact, we have $h^0(\omega_V(tL))>0$ for some integer $t$ with $0<t\le n$ for
arbitrary normal $n$-fold $V$, unless $V\cong\BP^n$.

Indeed, if not, we have $0=h^0(M,K_M+tL)=h^n(M,-tL_M)$ for any smooth model
$M\lra V$, so [{\bf F1};(2.2)] applies.

\Rmk By the same result, we can prove $\ke(V,B,L)\ge-n-1$ for arbitrary nef big
line bundle $L$.

\enddocument